\documentclass[11pt,a4paper]{article}
\usepackage{latexsym}
\usepackage{graphicx}
\usepackage{amsmath}
\usepackage{amssymb}
\usepackage{epsfig}

\makeatletter \@addtoreset{equation}{section}
 \makeatother

\begin{document}

\title{ \textbf{ On $\gamma$-Semi-Continuous Functions}}
\author{Sabir $Hussain^{1}$, Bashir $Ahmad^{2}$ and Takashi $Noiri^{3}$ \\
$^{1}$Department of Mathematics, Islamia University Bahawalpur,
Pakistan.\\ {\bf Present Address:} Department of Mathematics,
Yanbu University,\\ P. O. Box 31387, Yanbu,  Saudi Arabia.\\ E.
mail: sabiriub@yahoo.com.
\\ $^{2}$Department of Mathematics,
King Abdul Aziz University, P. O. Box 80203,\\ Jeddah 21589, Saudi
Arabia.\\ E. mail: drbashir9@gmail.com.\\
$^{3}$2949-1 Shiokita-cha, Hinagu, Yatsushiro-shi, Kumamoto-ken,
869-5142, Japan.\\ E. mail: t.noiri@nifty.com.\\}
\date{}
\maketitle
\textbf{Abstract.} In this paper, we continue studying the
properties of   $\gamma$-semi-continuous and $\gamma$-semi-open
functions introduced in [5].
AMS(2000)Subject Classification. Primary 54A05, 54A10, 54D10.
\section {Introduction}
     N. Levine [11] introduced the notion of semi-open sets in
     topological spaces. A. Csaszar [7,8] defined generalized open sets in
     generalized
topological spaces. In 1975, Maheshwari and Prasad [12] introduced
concepts of semi-$T_1$-spaces and semi-$R_0$-spaces. In 1979, S.
Kasahara [10] defined an operation $\alpha$  on topological
spaces. Carpintero, et. al [6] introduced the notion of
$\alpha$-semi-open sets as a generalization of semi-open sets. B.
Ahmad and F.U. Rehman [1, 14] introduced the notions of
$\gamma$-interior, $\gamma$-boundary and $\gamma$-exterior points
in topological spaces. They also studied properties and
characterizations of ($\gamma$ ,$\beta$ )-continuous mappings
introduced by H. Ogata [13]. In [2-4], B. Ahmad and S. Hussain
introduced the concept of $\gamma_ 0$-compact, $\gamma^*$-regular,
$\gamma$-normal spaces and explored their many interesting
properties. They initiated and discussed the concepts of
$\gamma^*$-semi-open sets , $\gamma^*$-semi-closed sets,
$\gamma^*$-semi-closure, $\gamma^*$-semi-interior points in
topological spaces [5,9]. In [9], they  introduced
$\Lambda_{s}^{\gamma}$-set and $\Lambda^{s^{\gamma}}$-set by using
$\gamma^*$-semi-open sets. Moreover, they also introduced
$\gamma$-semi-continuous function and $\gamma$-semi-open (closed)
functions in topological spaces and established several
interesting properties.

      In this paper, we continue studying the properties of  $\gamma$-semi-continuous functions
 and  $\gamma$-semi-open function introduced by B. Ahmad and S. Hussain
 [5].

        Hereafter, we shall write space in place of topological space in
the sequel.\\
\section {Preliminaries}

    We recall some definitions and results used in this paper to
    make it self-contained.\\
{\bf Definition [13].}  Let (X,$\tau$) be a space. An operation
$\gamma$ : $\tau\rightarrow$ P(X) is a function from $\tau$ to the
power set of X such that V $\subseteq V^\gamma$ , for each V
$\in\tau$, where $V^\gamma$ denotes the value of $\gamma$ at V.
The operations defined by $\gamma$(G) = G, $\gamma$(G) = cl(G) and
$\gamma$(G) = intcl(G) are examples of operation $\gamma$.\\
{\bf Definition [13].}  Let A be a subset of a space X. A point x
$\in$ A is said to be $\gamma$-interior point of A, if there
exists an open nbd N of x such that $ N^\gamma\subseteq$ A. The
set of all such points is denoted by $int_\gamma$(A). Thus
\begin{center}
        $int_\gamma$ (A) = $\{ x \in A : x \in N \in \tau$  and  $N^\gamma\subseteq A \} \subseteq A$.
\end{center}

    Note that A is $\gamma$-open [13] iff A =$int_\gamma$(A). A set A is
called $\gamma$- closed [13] iff X-A is $\gamma$-open.\\
{\bf Definition [10].}  A point x $\in$ X is called a
$\gamma$-closure point of A $\subseteq$ X, if $U^\gamma\cap A \neq
\phi$, for each open nbd U of x. The set of all $\gamma$-closure
points of A is called $\gamma$-closure of A and is denoted by
$cl_\gamma$(A). A subset A of X is called $\gamma$-closed, if
$cl_\gamma(A)\subseteq$A. Note that $cl_\gamma$(A) is contained in
every $\gamma$-closed superset
of A.\\
{\bf Definition [14].} The $\gamma$-exterior of A, written
$ext_{\gamma}(A)$ is defined as the $\gamma$-interior of $(X -
A)$. That is, $int_{\gamma} (A) = ext_{\gamma}(X - A)$.\\
{\bf Definition [14].} The $\gamma$-boundary of A, written
$bd_{\gamma}(A)$ is defined as the set of points which do not
belong to $\gamma$-interior or the $\gamma$-exterior of A.\\
{\bf Definition [13].} An operation $\gamma$ on $\tau$ is said be
regular, if for any open nbds U,V of x $\in$ X, there exists an
open nbd W of x such that $U^\gamma \cap V^\gamma\supseteq
W^\gamma$.\\
{\bf Definition [13].}  An operation $\gamma$ on $\tau$ is said to
be open, if for every open nbd U of each $x \in X$, there exists
$\gamma$-open set B such that $x\in B$ and $U^\gamma\subseteq
B$.\\
{\bf Definition [2].} Let  $A \subseteq X$. A point $x \in X$  is
said to be $\gamma$-limit point of A, if  $U \cap \{A - \{x\} \}
\neq \phi$, where U is a $\gamma$-open set in X. The set of all
$\gamma$-limit points
of A denoted  $A^{d}_{\gamma}$ is called $\gamma$-derived set.\\
{\bf Definition [9].} A subset A of a space (X,$\tau$ ) is said to
be a $\gamma^{*}$-semi-open set, if  there exists a $\gamma$ -open
set O such that $O \subseteq A \subseteq cl_{\gamma} (O)$. The set
of all $\gamma^{*}$-semi-open sets is denoted by
$SO_{\gamma^{*}}(X)$.\\
{\bf Definition [5].} A function $f: (X, \tau) \rightarrow (Y,
\tau)$ is said to be $\gamma$-semi-continuous,
if for any  $\gamma$-open B of Y, $f^{-1}(B)$ is  $\gamma^{*}$-semi-open in X.\\
{\bf Definition [5].} A function $f: (X, \tau) \rightarrow (Y,
\tau)$  is said to be $\gamma$-semi-open (closed), if for each
$\gamma$-open (closed) set U in X, $f(U)$  is
$\gamma^{*}$-semi-open
(closed) in Y.\\
{\bf Definition [5].}  A set A in a space X is said to be
$\gamma^{*}$-semi-closed, if there exists a  $\gamma$-closed set F
such that $int_{\gamma}(F) \subseteq A \subseteq F$.\\
{\bf Proposition [5].}  A subset A of X is
$\gamma^{*}$-semi-closed if and only if $X -A$
is  $\gamma^{*}$-semi-open.\\
{\bf Definition 2.1.}  A subset A of X is said to be
$\gamma$-semi-nbd of a point $x \in X$,  if there exists a  $\gamma^{*}$-semi-open set U such that  $x \in U \subseteq A$ .\\
{\bf Definition [9].}  Let A be a subset of space X . The
intersection of  all  $\gamma^{*}$-semi-closed sets containing A
is called $\gamma^{*}$-semi-closure of A and is denoted by
$scl_{\gamma^{*}} (A)$. Note that A is $\gamma^{*}$-semi-closed if
and
only if  $scl_{\gamma^{*}} (A) = A$.\\
{\bf Definition [5].} Let A be a subset of a space X . The union
of all $\gamma^{*}$-semi-open sets of X contained in A is called
$\gamma^{*}$-semi-interior of
A and is denoted by $sint_{\gamma^{*}} (A)$ .\\
{\bf Lemma 2.2.}  Let A be a subset of a space X. Then $x \in
scl_{\gamma^{*}} (A)$  if and only if for any
 $\gamma$-semi-nbd $N_{x}$  of x in X, $A \cap N_{x} \neq \phi$.\\
{\bf Proof. } Let $x \in scl_{\gamma^{*}} (A)$ .Suppose on the
contrary, there exists a $\gamma$-semi-nbd  $N_{x}$ of x in X such
that $A \cap N_{x} \neq \phi$. Then there exists $U \in
SO_{\gamma^{*}} (A)$  such that $x \in U \subseteq N_{x}$.
Therefore, $U \cap A = \phi$, so that $A \subseteq X -U$. Clearly
$X - U$ is $\gamma^{*}$-semi-closed in X and hence
$scl_{\gamma^{*}} (A) \subseteq X - U$. Since  $x \notin X - U$,
we obtain $x \notin scl_{\gamma^{*}} (A)$. This is contradiction
to the hypothesis. This proves the necessity.

              Conversely, suppose that every  $\gamma$-semi-nbd of x in X meets A. If  $x \notin scl_{\gamma^{*}} (A)$, then by
definition there exists a  $\gamma^{*}$-semi-closed F of X such
that  $A \subseteq F$ and  $x \notin F$. Therefore we have $x \in
X -F \in SO_{\gamma^{*}} (X)$. Hence $X - F$ is $\gamma$-semi-nbd
of x in X. But $(X -F) \cap A = \phi$. This is contradiction to
the hypothesis. Thus $x \in scl_{\gamma^{*}} (A)$.

\section {$\gamma$-Semi-Open Functions}

{\bf Theorem 3.1.}  Let  $f: X \rightarrow Y$ be a function from a
space X into a space Y and $\gamma$ is an open, monotone and
regular operation. Then the
following statements are equivalent:\\
    (1)  $f$  is  $\gamma$-semi-open.\\
    (2)   $ f(int_{\gamma}(A)) \subseteq sin_{\gamma^{*}} (f((A))$ for each subset A of X.\\
    (3)   For each  $x \in X$ and each  $\gamma$-open-nbd U of x, there exists a  $\gamma$-semi-nbd V of
          $f(x)$  such that $V \subseteq f(U)$.\\
{\bf Proof. } $(1) \Rightarrow (2)$. Suppose that $f$  is
$\gamma$-semi-open, and let A be an arbitrary subset of X. Since $
f(int_{\gamma}(A))$ is $\gamma^{*}$-semi-open and $
f(int_{\gamma}(A)) \subseteq f(A)$, then $ f(int_{\gamma}(A))
\subseteq sin_{\gamma^{*}} (f((A))$.

        $(2) \Rightarrow (3)$. Let U be an arbitrary  $\gamma$-open-nbd of $x \in X$. Then there exists  $\gamma$-open set O such that
$x \in O \subseteq U$. By hypothesis, we have $ f(O) =
f(int_{\gamma}(O)) \subseteq sin_{\gamma^{*}} (f((O))$  and hence
$ f(O) \subseteq sin_{\gamma^{*}} (f((O))$. Therefore it follows
that $f(O)$  is $\gamma$-semi-open-nbd in Y such that $f(x) \in
f(O) \subseteq f(U)$. This proves (3).

        $(3) \Rightarrow (1)$. Let U be an arbitrary  $\gamma$-open set in X. For each $y \in f(U)$, by hypothesis
there exists a  $\gamma$-semi-nbd   $V_{y}$ of y in Y such that
$V_{y} \subseteq f(U)$. Since $V_{y}$ is a $\gamma$-semi-nbd of y,
there exists a $\gamma^{*}$-semi-open set $A_{y}$ in Y such that
$y \in A_{y} \subseteq V_{y}$. Therefore  $f(U) = \bigcup \{A_{y}
: y \in f(U) \}$ is a $\gamma^{*}$-semi-open in Y, since is
$\gamma$ regular [9]. This
shows that  $f$ is a  $\gamma$-semi-open function.\\
\\
{\bf Theorem 3.2.} A bijective function $f: X \rightarrow Y$   is
$\gamma$-semi-open if and only if $ f^{-1} (scl_{\gamma^{*}}(B))
\subseteq cl_{\gamma}(f^{-1}(B))$ for every subset B of  Y,
where $\gamma$  is an open operation.\\
{\bf Proof.} Let B be an arbitrary subset of Y. Put
\begin{center}
 $U = X - cl_{\gamma} (f^{-1}(B))$  ...... (I)
\end{center}
Clearly U is a $\gamma$-open set in X. Then by hypothesis, $f(U)$
is a $\gamma^{*}$-semi-open set in Y, or $Y - f(U)$ is
$\gamma^{*}$-semi-closed set in Y. Since  $f$ is onto, from (I),
it follows   $B \subseteq Y - f(U)$. Thus we have
$scl_{\gamma^{*}}(B) \subseteq Y - f(U)$. Since $f$ is one-one, we
have $ f^{-1} (scl_{\gamma^{*}}(B)) \subseteq f^{-1}(Y) -
f^{-1}f(U) = X - f^{-1}f(U) \subseteq X - U = cl_{\gamma} (f^{-1}
(B))$. This proves the necessity.

       Conversely, let U be an arbitrary  -open set in X. Put  $B = Y - f(U)$. Since  $f$ is bijective, therefore
by hypothesis,  $f(U) \cap scl_{\gamma^{*}}(B) = f(U \cap
f^{-1}(scl_{\gamma^{*}}(B))) \subseteq f(U \cap
cl_{\gamma}(f^{-1}(B)))$. Since U is $\gamma$-open, therefore by
Lemma 2(3) [14], we have $U \cap cl_{\gamma}(f^{-1}(B)) \subseteq
cl_{\gamma}(U \cap f^{-1}(B))$. Moreover, it is obvious that $U
\cap f^{-1}(B) = \phi$. Thus we have $f(U) \cap
scl_{\gamma^{*}}(B) = \phi$ and hence $scl_{\gamma^{*}}(B)
\subseteq Y - f(U) = B$. Therefore B is a $\gamma^{*}$-semi-closed
in Y and hence $f(U)$ is a $\gamma^{*}$-semi-open set in Y. This
proves that f is a $\gamma$-semi-open
mapping.\\
\\
{\bf Definition 3.3 [13].} A function  $f: (X, \tau, \gamma)
\rightarrow (Y, \delta, \beta)$ is said to be $(\gamma ,
\beta)$-continuous, if for each $x \in X$ and each open set V
containing $f(x)$, there exists an open set U such that $x \in U$
and $f(U^{\gamma}) \subseteq V^{\beta}$, where $\gamma$ and
$\beta$ are operations on $\tau$ and $\delta$
respectively.\\
\\
{\bf Definition 3.4 [13].} A function  $f: (X, \tau, \gamma)
\rightarrow (Y, \delta, \beta)$ is said to be $(\gamma ,
\beta)$-open (closed), if for any $\gamma$-open (closed) set A of
X, $f(A)$ is $\gamma$-open
(closed) in Y.\\
\\
{\bf Theorem 3.5 [1].}  Let  $f: (X, \tau, \gamma) \rightarrow (Y,
\delta, \beta)$ be a function and  $\beta$ be an open operation on
Y. Then $f$  is  $(\gamma , \beta)$-continuous if and only if for
each
$\beta$-open set V in Y, $f^{-1}(V)$  is  $\gamma$-open in X.\\
\\
{\bf Theorem 3.6 [1].} Let $f: (X, \tau, \gamma) \rightarrow (Y,
\delta, \beta)$  be a function and $\beta$  be an open
operation on Y. Then the following are equivalent:\\
    (1) $f$ is  $(\gamma , \beta)$-open .\\
    (2) $f^{-1}(cl_{\beta} (B)) \subseteq cl_{\gamma}(f^{-1}(B))$.\\
    (3) $f^{-1}(bd_{\beta} (B)) \subseteq bd_{\gamma}(f^{-1}(B))$ for any subset B of Y.\\
\\
{\bf Theorem 3.7.}  If a function  $f: (X, \tau, \gamma)
\rightarrow (Y, \delta, \beta)$ is a  $(\gamma , \beta)$-open and
a $(\gamma , \beta)$-continuous, then the inverse image
$f^{-1}(B)$ of each $\beta^{*}$-semi-open set B in Y is
$\gamma^{*}$-semi-open in X, where $\beta$  is an open operation on Y.\\
{\bf Proof.}  Let B be an arbitrary  $\beta^{*}$-semi-open set in
Y. Then there exists  $\beta$-open set V in Y such that  $V
\subseteq B \subseteq cl_{\beta} (V)$. Since f is $(\gamma ,
\beta)$-open, using Theorem 3.6, we have  $f^{-1}(V) \subseteq
f^{-1}(B) \subseteq f^{-1}(cl_{\beta} (V)) \subseteq
cl_{\gamma}(f^{-1}(V))$. Since is $(\gamma , \beta)$-continuous
and V is $\beta$-open in Y, by Theorem 3.5, $f^{-1}(V)$ is
$\gamma$-open in X. This shows that $f^{-1}(B)$ is
$\gamma^{*}$-semi-open set in X.\\
\\
{\bf Theorem 3.8.}  Let X, Y and Z be three spaces and let $f: X
\rightarrow Y$ be a function, $g: Y \rightarrow Z$  be an
injective function and  $gof: X \rightarrow Z$ is a
$\gamma$-semi-open
function. Then we have:\\
(1) If f is  $(\gamma , \beta)$-continuous and surjective, then g is  $\gamma$-semi-open.\\
(2) If g is  $(\beta, \alpha)$-open, $(\beta, \alpha)$-continuous
and injective, then f is $\gamma$-semi-open, where $\beta$  is open operation on Y.\\
{\bf Proof.}  (1)  Let V be a  $\beta$-open set in Y. Then
$f^{-1}(V)$  is $\gamma$-open in X, because $f$ is  $(\gamma ,
\beta)$-continuous. Since $gof$  is $\gamma$-semi-open and $f$ is
surjective, therefore $g(V) = (gof)(f^{-1}(V))$ is
$\alpha^{*}$-semi-open in Z. This shows that g is
a  $\gamma$-semi-open function.\\
 (2)   Since g is injective, therefore
for $A \subseteq X$,  $f(A) = g^{-1}(g(f(A)))$. Let U be a
 $\gamma$-open set in X, then  $gof(U)$  is  $\alpha^{*}$-semi-open. Thus
 by Theorem 3.7,  $g^{-1}(g(f(U))) = f(U)$ is  $\beta^{*}$-semi-open in Y. This shows that $f$ is a  $\gamma$-semi-open
 function.\\

    Let  $B \subseteq X$, $\gamma : \tau
\rightarrow P(X)$  be an operation. We define  $\gamma_{B} :
\tau_{B} \rightarrow P(X)$  as  $\gamma_{B}(U \cap B) = \gamma(U)
\cap B$. From here $\gamma_{B}$ is an operation and satisfies that
$cl_{\gamma_B}(U \cap B) \subseteq cl_{\gamma}(U \cap B) \subseteq
cl_{\gamma}(U) \cap cl_{\gamma}(B)$. Using this fact we prove the following:\\
\\
{\bf Theorem 3.9.}  Let X be a space and B a
$\gamma^{*}$-semi-open set in X containing a subset A of X. If A
is  $\gamma^{*}$-semi-open in the subspace
B, then A is  $\gamma^{*}$-semi-open in X, where  $\gamma$ is a regular operation.\\
{\bf Proof.}  Let A be  $\gamma^{*}_{B}$-semi-open in the subspace
B. Then there exists a  $\gamma_{B}$-open set $U_{B}$  in B such
that $U_{B} \subseteq A \subseteq cl_{\gamma_{B}}(U_{B})$. Since
$U_{B}$ is $\gamma_{B}$-open in B, there exists a $\gamma$-open
set U in X such that $U_{B} = U \cap B$[4]. Thus we have $U \cap B
\subseteq A \subseteq cl_{\gamma_{B}} (U \cap B) \subseteq
cl_{\gamma}(U \cap B) = cl_{\gamma}(A) \cap cl_{\gamma}(B)$. Since
B is $\gamma^{*}$-semi-open set in X and U is $\gamma$-open in X,
therefore $U \cap B$
is $\gamma$-open in X. Consequently, A is a  $\gamma^{*}$-semi-open set in X.\\
\\
{\bf Theorem 3.10.}  Let X and Y be spaces. If a bijective
function $f: X \rightarrow Y$ is a  $\gamma$-semi-open, then for
each  $\gamma$-open set $V(\neq \phi )$ in Y $f |_{f^{-1}(V)}:
f^{-1}(V) \rightarrow V$ is
$\gamma$-semi-open, where $\gamma$  is a regular operation.\\
{\bf  Proof.}  Let   $U_{V}$ be an arbitrary
$\gamma_{f^{-1}(V)}$-open set in $f^{-1}(V)$. Then there exists a
$\gamma$-open set U in X such that $U_{V} = U \cap f^{-1}(V)$. Now
we have $[f |_{f^{-1}(V)}] (U_{V}) =f(U \cap f^{-1}(V)) = f(U)
\cap V$. Since $f(U)$ is $\gamma^{*}$-semi-open and V is
 $\gamma$-open,  $f(U) \cap V$ is  $\gamma^{*}$-semi-open. Hence $[f |_{f^{-1}(V)}](U_{V})$ is also  $\gamma^{*}_{V}$-semi-open in V.
 This shows that $f |_{f^{-1}(V)}:
f^{-1}(V) \rightarrow V$  is a  $\gamma$-semi-open
 mapping.\\
\\
{\bf Theorem 3.11.} A bijective function  $f: X \rightarrow Y$ is
$\gamma$-semi-open if and only if for any subset V of Y and for
any  $\gamma$-closed set F of X containing  $f^{-1}(V)$, there
exists a $\gamma^{*}$-semi-closed set G of Y containing
V such that  $f^{-1}(G) \subseteq F$.\\
{\bf Proof.} Let $V \subseteq Y$ and F be a  $\gamma$-closed set
of X containing $f^{-1}(V)$. Put  $G = Y - f(X - F)$. Since  $f$
is $\gamma$-semi-open, so G is $\gamma^{*}$-semi-closed sets in Y
. As $f$  is bijective, it follows from  $f^{-1}(V) \subseteq F$
that $V \subseteq G$. Calculations give $f^{-1}(G) \subseteq F$.

       Conversely, suppose U is  $\gamma$-open set. Put $V = Y - f(U)$. Then  $X - U$ is  $\gamma$-closed set in X
containing $f^{-1}(V)$. By hypothesis, there exists a
$\gamma^{*}$-semi-closed set G of Y such that $V \subseteq G$  and
$f^{-1}(G) \subseteq (X - U)$. On the other hand, it follows from
$V \subseteq G$ that $f(U) = (Y - V) \subseteq (Y - G)$.
Therefore, we obtain $f(U) = (Y - G) \in SO_{\gamma^{*}} (Y)$.
This shows
that $f$  is  $\gamma$-semi-open.\\
\\
{\bf Lemma 3.12 [5].}   The following properties of a subset A of
X are equivalent:\\
          (1) A is  $\gamma^{*}$-semi-closed.\\
          (2)  $int_{\gamma}(cl_{\gamma}(A)) \subseteq A$.\\
          (3)  $X - A$ is  $\gamma^{*}$-semi-open.\\
\\
{\bf Theorem 3.13.}  If  $f: X \rightarrow Y$ is  $(\gamma ,
\beta)$-open and  $(\gamma , \beta)$-continuous mapping, then the
inverse image $f^{-1}(B)$  of each $\beta^{*}$-semi-closed B in Y
is $\gamma^{*}$-semi-closed
in X, where  $\beta$ is an open operation on Y.\\
{\bf Proof.}  This follows from Theorem 3.7 and Lemma 3.12.\\
\\
{\bf Theorem 3.14.}  Let $f: X \rightarrow Y$  be surjective and
$g: Y \rightarrow Z$ be an injective function and let $gof: X
\rightarrow Z$ be a  $\gamma$-semi-closed function. Then \\
(1) If f is $(\gamma , \beta)$-continuous and surjective, then g is  $\beta$-semi-closed.\\
(2) If g is $( \beta, \alpha)$-open, $(\beta, \alpha)$-continuous
and injective, then $f$ is $\gamma$-semi-closed, where $\beta$  is an open operation on Y.\\
{\bf  Proof.} (1) Suppose H is an arbitrary $\beta$-closed set in
Y. Then $f^{-1}(H)$ is $\gamma$-closed in X because $f$ is
$(\gamma , \beta)$-continuous. Since $gof$ is $\gamma$-semi-closed
and $f$ is surjective,  $gof(f^{-1}(H)) \subseteq g(f(f^{-1}(H)))
= g(H)$, is $\alpha^{*}$-semi-closed in Z. This implies that g is
 $\beta$-semi-open function. This proves (1).\\
         (2) Since g is injective so for every subset A of X,   $f(A) = g^{-1}(g(f(A)))$. Let F be an
 arbitrary  $\gamma$-closed set in X. Then $gof(F)$  is  $\gamma^{*}$-semi-closed. It follows immediately from
 Theorem 3.13 that $f(F)$ is  $\gamma^{*}$-semi-closed set in Y. This implies that f is  $\gamma$-semi-closed.

\section {$\gamma$-Semi-Closed Functions}
{\bf Theorem 4.1.} Let $\gamma$  be an open and monotone
operation. A function $f: X \rightarrow Y$ is $\gamma$-semi-closed
if and only if $f(cl_{\gamma}(A)) \supseteq
int_{\gamma}(cl_{\gamma}(f(A)))$ for every subset A of
X. \\
{\bf Proof.} Suppose $f$ is a  $\gamma$-semi-closed mapping and A
is an arbitrary subset of X. Then $f(cl_{\gamma}(A))$  is
$\gamma^{*}$-semi-closed in Y. Then by Lemma 3.12, we obtain
$f(cl_{\gamma}(A)) \supseteq
int_{\gamma}(cl_{\gamma}(f(cl_{\gamma}(A)))) \supseteq
int_{\gamma}(cl_{\gamma}(f(A)))$. This implies that
$f(cl_{\gamma}(A)) \supseteq int_{\gamma}(cl_{\gamma}(f(A)))$.

       Conversely, suppose that F is an arbitrary  $\gamma$-closed set in X. Then by hypothesis, we
have $int_{\gamma}(cl_{\gamma}(f(F))) \subseteq f(cl_{\gamma}(F))
= f(F)$. By Lemma 3.12, $f(F)$ is $\gamma^{*}$-semi-closed in Y.
This implies that
$f$ is  $\gamma$-semi-closed.\\

    Recall [9] that the intersection of all  $\gamma^{*}$-semi-closed sets
containing A is called
 $\gamma$-semi-closure of A and is denoted by  $scl_{\gamma^{*}}(A)$. Clearly A is  $\gamma^{*}$-semi-closed if and only if
 $scl_{\gamma^{*}}(A) = A$.\\
\\
{\bf Theorem 4.2.}  Let  $\gamma$   be an open and monotone
operation. A function  $f: X \rightarrow Y$ is
$\gamma$-semi-closed if and only if $scl_{\gamma^{*}}(A) \subseteq
f(cl_{\gamma}(A))$ for every subset A of
X.\\
{\bf Proof.} Suppose  $f$ is a  $\gamma$-semi-closed mapping and A
is an arbitrary subset of X. Then $f(cl_{\gamma}(A))$  is
$\gamma^{*}$-semi-closed. Since $f(A) \subseteq
f(cl_{\gamma}(A))$, we obtain $scl_{\gamma^{*}}(f(A)) \subseteq
f(cl_{\gamma}(A))$. This implies $scl_{\gamma^{*}}(f(A)) \subseteq
f(cl_{\gamma}(A))$.

         Sufficiency follows from Theorem 4.1.\\
\\
{\bf Theorem 4.3.}  A surjective function $f: X \rightarrow Y$ is
$\gamma$-semi-closed if and only if for each subset B in Y and
each  $\gamma$-open set U in X containing  $f^{-1}(B)$, there
exists a $\gamma^{*}$-semi-open set V in Y containing B
such that  $f^{-1}(V) \subseteq U$, where $\gamma$  is a monotone and regular operation.\\
{\bf Proof.} Suppose B is an arbitrary subset in Y and U is an
arbitrary  $\gamma$-open set in X containing $f^{-1}(B)$. We put
\begin{center}
$V= Y - f(X - U)$ ..... (*)
\end{center}
Then V is $\gamma^{*}$-semi-open set in Y. Since $f^{-1}(B)
\subseteq U$, calculations give $B \subseteq V$. Moreover, by (*),
we have $f^{-1}(V) = f^{-1}(Y) - f^{-1}(f(X -U)) = X - f^{-1}(f(X
-U)) \subseteq X - (X - U) = U$.

    Conversely, suppose that F is an arbitrary  $\gamma$-closed set in X. Let y be an
arbitrary point in  $Y - f(F)$, then  $f^{-1}(y) \subseteq X -
f^{-1}(f(F)) \subseteq X - F$, and $X - F$ is $\gamma$-open in X.
Hence by the hypothesis, there exists a $\gamma^{*}$-semi-open set
$V_{y}$ containing y such that $f^{-1}(V_{y}) \subseteq X - F$.
This implies that $v \in V_{y} \subseteq Y - f(F)$. We obtain that
$Y - f(F) = \bigcup \{ V_{y} : y \in Y - f(F) \}$ is
$\gamma^{*}$-semi-open in Y, since union of any collection of
$\gamma^{*}$-semi-open sets is $\gamma^{*}$-semi-open. Therefore
$f(F)$ is $\gamma^{*}$-semi-closed.
\\
\section {$\gamma$-Semi-Continuous Functions}
{\bf Theorem 5.1.}  Let $f: X \rightarrow Y$  be a function and
$\gamma$  is an open operation. Then the following are equivalent:\\
 (1)    f is $\gamma$-semi-continuous.\\
 (2)  $int_{\gamma}(cl_{\gamma}(f^{-1}(B))) \subseteq f^{-1}(cl_{\gamma}(B))$ for each subset B of Y.\\
 (3)  $f(int_{\gamma}(cl_{\gamma}(A))) \subseteq cl_{\gamma}(f(A))$ for each subset A of X.\\
{\bf  Proof.} $(1) \Rightarrow (2)$. Let B be an arbitrary subset
of Y. Then by (1), $f^{-1}(cl_{\gamma}(B))$ is a
$\gamma^{*}$-semi-closed set of X. Since $B \subseteq
cl_{\gamma}(B)$, by Lemma 3.12, we get $f^{-1}(cl_{\gamma}(B))
\supseteq int_{\gamma}(cl_{\gamma}(f^{-1}(cl_{\gamma}(B))))
\supseteq int_{\gamma}(cl_{\gamma}(f^{-1}(B)))$implies that
$int_{\gamma}(cl_{\gamma}(f^{-1}(B))) \subseteq
f^{-1}(cl_{\gamma}(B))$.

            $(2) \Rightarrow (3)$. Let A be an arbitrary subset of X. Put $B = f(A)$. Then $A \subseteq f^{-1}(B)$.
Therefore by hypothesis, we have
$int_{\gamma}(cl_{\gamma}(A))\subseteq
int_{\gamma}(cl_{\gamma}(f^{-1}(B)))  \subseteq
f^{-1}(cl_{\gamma}(B))$. Consequently, we have
$f(int_{\gamma}(cl_{\gamma}(A)))\subseteq ff^{-1}(cl_{\gamma}(B))
\subseteq cl_{\gamma}(B) = cl_{\gamma} (f(A))$. This gives (3).

            $(3) \Rightarrow (1)$.  Let F be an arbitrary  $\gamma$-closed set of Y. Put $A = f^{-1}(F)$,
then  $f(A) \subseteq F$. Therefore by hypothesis, we have
\begin{center}
$f(int_{\gamma}(cl_{\gamma}(A))) \subseteq cl_{\gamma}(f(A))
\subseteq cl_{\gamma}(F) = F$ ..... (**)
\end{center}
By (**), we have $int_{\gamma}(cl_{\gamma}(A)) \subseteq
f^{-1}f(int_{\gamma}(cl_{\gamma}(A))) \subseteq
f^{-1}(cl_{\gamma}(f(A))) \subseteq f^{-1}(cl_{\gamma}(F)) =
f^{-1}(F)$, or $int_{\gamma}(cl_{\gamma}(A)) \subseteq
f^{-1}(F)$,. By Lemma 3.12, $f^{-1}(F)$ is a
$\gamma^{*}$-semi-closed set
in X. This implies that $f$ is  $\gamma$-semi-continuous.\\
\\
{\bf Definition 5.2.}   Let X be a space $A \subseteq X$ and  $p
\in X$. Then p is a $\gamma^{*}$-semi-limit point of A, for all
$\gamma^{*}$-semi-open set U containing p such that $U \cap (A -
\{p\}) \neq \phi$. The set of all $\gamma^{*}$-semi-limit point of
A is said to be $\gamma^{*}$-semi-derived set of A and is denoted
by $sd_{\gamma^{*}} (A)$.
\begin{center}
Clearly if $A \subseteq B$ then $sd_{\gamma^{*}} (A) \subseteq
sd_{\gamma^{*}} (B)$..... (I)
\end{center}
{\bf Remark 5.3.} From the definition, it follows that p is a
$\gamma^{*}$-semi-limit point of  A if and only if  $ p \in scl_{\gamma^{*}} (A - \{p\})$.\\
\\
{\bf Theorem 5.4.}  The  $\gamma^{*}$-semi-derived set,
$sd_{\gamma^{*}}$, has the
following properties:\\
    (1)  $scl_{\gamma^{*}} (A) = A \cup sd_{\gamma^{*}} (A)$.\\
    (2)  $sd_{\gamma^{*}} (A \cup B) = sd_{\gamma^{*}} (A) \cup sd_{\gamma^{*}} (B)$. In general\\
    (3)  $\bigcup_{i}sd_{\gamma^{*}} (A_{i}) = sd_{\gamma^{*}} (\bigcup_{i}(A_{i}))$.\\
    (4)   $sd_{\gamma^{*}} (sd_{\gamma^{*}} (A))\subseteq sd_{\gamma^{*}} (A)$.\\
    (5)  $scl_{\gamma^{*}} (sd_{\gamma^{*}} (A)) = sd_{\gamma^{*}} (A)$.\\
{\bf Proof.}   (1)  Let  $x \in scl_{\gamma^{*}} (A)$. Then $x \in
C$, for every $\gamma^{*}$-semi-closed
superset C of A. Now\\
    (i)  If $x \in A$, then $x \in A \cup sd_{\gamma^{*}} (A)$.\\
(ii)  If $x \notin A$, then we prove that $ x \in scl_{\gamma^{*}}
(A)$.\\
To prove (ii), suppose U is  $\gamma^{*}$-semi-open set containing
x. Then $U \cap A \neq \phi$, for otherwise, $A \subseteq X - U =
C$, where C is a $\gamma^{*}$-semi-closed superset of A not
containing x. This contradicts the fact that x belongs to every
$\gamma^{*}$-semi-closed superset C of A. Therefore $x \in
sd_{\gamma^{*}} (A)$ gives $x \in A \cup sd_{\gamma^{*}} (A)$.

    Conversely, suppose that  $x \in A \cup sd_{\gamma^{*}} (A)$,  we show that  $x \in scl_{\gamma^{*}} (A)$. If
    $x \in A$ then $x \in scl_{\gamma^{*}} (A)$. If $x \in sd_{\gamma^{*}} (A)$, then we show that x  is in every
$\gamma^{*}$-semi-closed superset of A. We suppose otherwise that
there is $\gamma^{*}$-semi-closed superset C  of A not containing
x. Then $x \in X- C = U$(say), which is $\gamma^{*}$-semi-open and
$U \cap A = \phi$.  This implies that $x \notin sd_{\gamma^{*}}
(A)$. This contradiction proves that $x \in scl_{\gamma^{*}} (A)$.
Consequently $scl_{\gamma^{*}} (A) = A \cup sd_{\gamma^{*}} (A)$.
This proves (1).

       (2) $sd_{\gamma^{*}} (A \cup B) \subseteq sd_{\gamma^{*}} (A) \cup sd_{\gamma^{*}} (B)$.\\
Let $x \in sd_{\gamma^{*}} (A \cup B)$. Then $x \in
scl_{\gamma^{*}} ((A \cup B) - \{x\})$ or
 $x \in
scl_{\gamma^{*}} ((A - \{x\}) \cup (B - \{x\})$ implies  $x \in
scl_{\gamma^{*}} (A  - \{x\})$ or $x \in scl_{\gamma^{*}} (B -
\{x\})$. This gives $x \in sd_{\gamma^{*}} (A)$ or $x \in
sd_{\gamma^{*}} (B)$ . Therefore $x \in sd_{\gamma^{*}} (A)\cup
sd_{\gamma^{*}} (B)$ . This
  proves  $sd_{\gamma^{*}} (A \cup B) \subseteq sd_{\gamma^{*}} (A) \cup sd_{\gamma^{*}} (B)$

             Converse follows directly by using the property(I).

       (3)  The proof is immediate by property (I).

       (4)  Suppose that $x \notin  sd_{\gamma^{*}} (A)$. Then $x \notin  scl_{\gamma^{*}} (A - \{x\})$. This implies that there is
 $\gamma^{*}$-semi-open set U such that $x \in U$ and $U \cap (A - \{x \}) = \phi$. We prove that
$x \notin  sd_{\gamma^{*}} (sd_{\gamma^{*}} (A))$. Suppose on the
contrary that $x \in  sd_{\gamma^{*}} (sd_{\gamma^{*}} (A))$. Then
$x \in  scl_{\gamma^{*}} (sd_{\gamma^{*}} (A) - \{x \})$. Since $x
\in U$, we have $U \cap (sd_{\gamma^{*}} (A) - \{x \}) \neq \phi$.
Therefore there is a $q \neq x$ such that $q \in  U  \cap
(sd_{\gamma^{*}} (A))$. It follows that $q \in (U - \{x\}) \cap
(sd_{\gamma^{*}} (A) -\{x\})$. Hence $(U - \{x\}) \cap (
sd_{\gamma^{*}} (A) -\{x\}) \neq \phi$, a contradiction to the
fact that $(U \cap ( sd_{\gamma^{*}} (A) -\{x\}) = \phi$. This
implies that
 $x \notin sd_{\gamma^{*}} (sd_{\gamma^{*}} (A))$  and so  $sd_{\gamma^{*}} (sd_{\gamma^{*}} (A))\subseteq sd_{\gamma^{*}} (A)$.
 This proves (4).

       (5)   This is a consequence of (1), (2) and (4).\\
\\
{\bf Theorem 5.5 [5].}  Let $f: X \rightarrow Y$  be a function.
Then the following
are equivalent:\\
   (1)   $f: X \rightarrow Y$   is  $\gamma$-semi-continuous.\\
   (2)    $scl_{\gamma^{*}} (f^{-1}(A)) \subseteq f^{-1}(cl_{\gamma}(A))$  for each subset A of Y.\\
\\
{\bf Theorem 5.6.} Let $f: X \rightarrow Y$  be a function and
$\gamma$ is an open
operation. Then the following are equivalent:\\
   (1)  $f: X \rightarrow Y$    is  $\gamma$-semi-continuous.\\
   (2)  $f(sd_{\gamma^{*}} (A)) \subseteq cl_{\gamma}(f(A))$ for any subset A of X.\\
{\bf Proof.}  $(1) \Rightarrow (2)$. Suppose that $f$  is
$\gamma$-semi-continuous. Let A be any set in X. Since
$cl_{\gamma}(f(A))$ is $\gamma$-closed in Y.
$f^{-1}(cl_{\gamma}(A))$ is $\gamma^{*}$-semi-closed in X. $A
\subseteq f^{-1}(f(A)) \subseteq f^{-1}(cl_{\gamma}(f(A)))$ gives
$scl_{\gamma^{*}}(A) \subseteq
scl_{\gamma^{*}}(f^{-1}(cl_{\gamma}(f(A)))) =
f^{-1}(cl_{\gamma}(f(A)))$. Therefore $f(sd_{\gamma^{*}}(A))
\subseteq f(scl_{\gamma^{*}} (A)) \subseteq f
f^{-1}(cl_{\gamma}(f(A))) \subseteq cl_{\gamma}(f(A))$.
Consequently, $f(sd_{\gamma^{*}} (A)) \subseteq
cl_{\gamma}(f(A))$.

     $(2) \Rightarrow (1)$.  Suppose that $f(sd_{\gamma^{*}} (A)) \subseteq cl_{\gamma}(f(A))$, for $A \subseteq X$. Let B be any
 $\gamma$-closed subset of Y. We show that $f^{-1}(B)$  is  $\gamma^{*}$-semi-closed in X. By hypothesis,
 $f(sd_{\gamma^{*}} (f^{-1}(B))) \subseteq cl_{\gamma}(f(f^{-1}(B))) \subseteq cl_{\gamma}(B) = B$  or
 $f(sd_{\gamma^{*}} (f^{-1}(B))) \subseteq  B$  gives
 $sd_{\gamma^{*}} (f^{-1}(B)) \subseteq f^{-1}f(sd_{\gamma^{*}} (f^{-1}(B))) \subseteq f^{-1}(B)$ or
 $sd_{\gamma^{*}} (f^{-1}(B)) \subseteq  f^{-1}(B)$ implies $f^{-1}(B)$  is  $\gamma^{*}$-semi-closed in X. Thus $f$ is  $\gamma$-semi-continuous.\\
\\
{\bf Theorem 5.7 [5].} Let $f: X \rightarrow Y$ be a function and
$x \in X$. Then $f$ is $\gamma$-semi-continuous if and only if for
each $\gamma$-open set B containing f(x), there exists $A \in
SO_{\gamma^{*}}(X)$ such that $x \in A$ and $f(A) \subseteq B$,
where $\gamma$ is a regular operation.\\

    We use Theorem 5.7 and prove the following:\\
{\bf Theorem 5.8.}  Let $f: X \rightarrow Y$  be an injective
function. If $f$ is $\gamma$-semi-continuous then
$f(sd_{\gamma^{*}}(A)) \subseteq (f(A))^{d}_{\gamma}$ for every
$A \subseteq X$, where $\gamma$ is a regular
operation.\\
{\bf Proof.} Suppose that $f$ is  $\gamma$-semi-continuous. Let $A
\subseteq X$, $a \in sd_{\gamma^{*}}(X)$ and V be a
$\gamma$-open-nbd of $f(a)$. Since $f$ is $\gamma$-semi-continuous
then by Theorem 5.7, there exists a
 $\gamma$-semi-open-nbd U of a such that  $f(U) \subseteq V$. But  $a \in sd_{\gamma^{*}}(A)$, therefore there exists an element
 $a_{1} \in U \cap A$  such that  $a \neq a_{1}$; then $f(a_{1}) \in f(A)$  and since f is an injection $f(a) \neq f(a_{1})$. Thus every
$\gamma$-open-nbd V of $f(a)$ contains an
 element  $f(a_{1})$ of $f(A)$ different from  $f(a)$. Consequently  $f(a) \in (f(A))^{d}_{\gamma}$. We have therefore,
 $f(sd_{\gamma^{*}}(A)) \subseteq (f(A))^{d}_{\gamma}$.\\

         The following theorem follows from Theorem 5.6:\\
{\bf Theorem 5.9.}  Let $f: X \rightarrow Y$  be a function. If
for every  $A \subseteq X$,  $f(sd_{\gamma^{*}}(A)) \subseteq
(f(A))^{d}_{\gamma}$, then $f$ is $\gamma$-semi-continuous, where
$\gamma$ is an open operation.
\\
{\bf Theorem 5.10.}  A function $f: X \rightarrow Y$  is
$\gamma$-semi-continuous if
and only if $f^{-1}(int_{\gamma}(B)) \subseteq sint_{\gamma^{*}}(f^{-1}(B))))$, for each $B \subseteq Y$, where  $\gamma$ is a regular operation.\\
{\bf Proof.} For any $B \subseteq Y$, $int_{\gamma}(B) = Y -
cl_{\gamma}(Y - B)$ [14]. This implies  $f^{-1}(int_{\gamma}(B)) =
f^{-1}(Y - cl_{\gamma}(Y - B)) = X - f^{-1}(cl_{\gamma}(Y - B))$.
Since f is $\gamma$-semi-continuous, by Theorem 5.5 we have
$scl_{\gamma^{*}}(f^{-1}(Y - B)) \subseteq f^{-1}(cl_{\gamma}(Y
-B))$. Hence  $f^{-1}(int_{\gamma}(B)) \subseteq X -
scl_{\gamma^{*}}(f^{-1}(Y - B))$. Thus $f^{-1}(int_{\gamma}(B))
\subseteq X - scl_{\gamma^{*}}(X - f^{-1}(B))$. Hence
$f^{-1}(int_{\gamma}(B)) \subseteq X - scl_{\gamma^{*}}(X -
f^{-1}(B)) = sint_{\gamma^{*}}(f^{-1}(B))))$.

     Conversely, let B be an arbitrary  $\gamma$-open set in Y. Then $int_{\gamma}(B) = B$. By
hypothesis  $f^{-1}(B) = f^{-1}(int_{\gamma}(B)) \subseteq
sint_{\gamma^{*}}(f^{-1}(B))$ implies $f^{-1}(B)  \subseteq
sint_{\gamma^{*}}(f^{-1}(B))$. But $sint_{\gamma^{*}}(f^{-1}(B))
\subseteq f^{-1}(B)$. Therefore,  $f^{-1}(B) =
sint_{\gamma^{*}}(f^{-1}(B))$. Thus $f^{-1}(B)$
is  $\gamma^{*}$-semi-open . Consequently, $f$  is  $\gamma$-semi-continuous.\\

\end{document}